\documentclass[a4paper,11pt,reqno]{amsart}

\usepackage[utf8]{inputenc}		% font and encoding
\usepackage[T1]{fontenc}
\usepackage[english]{babel}

\usepackage{amsfonts}			% ams-packages
\usepackage{amsmath}
\usepackage{amssymb}
\usepackage{amsthm}

\usepackage{mathtools}

\usepackage[normalem]{ulem}

\usepackage{graphicx}	
\usepackage{tikz}
\usetikzlibrary{arrows.meta}
\usetikzlibrary
{arrows.meta,decorations.pathreplacing}
\colorlet{leftpart}{blue}
\colorlet{middlepart}{red}
\colorlet{rightpart}{blue}

\definecolor{curve}{RGB}{255,0,0}
\definecolor{endpoint}{RGB}{0,0,255}
\colorlet{annulusfill}{black!10}
\colorlet{boundary}{black}
\colorlet{guide}{black}% images and drawing
\usetikzlibrary{cd}
\usetikzlibrary{arrows.meta}

\usepackage{hyperref}			% hyperlinks
\hypersetup{colorlinks=false, urlcolor=black, linkcolor=black}
\usepackage{cleveref}

\usepackage{enumitem}			% arbitrary list items

\usepackage{color}				% colored text for signaling changes
\usepackage{float}

\usepackage{mathrsfs}

\newcommand{\N}{\mathbb{N}}						% Nonnegative integers
\newcommand{\R}{\mathbb{R}}						% Reals
\newcommand{\C}{\mathbb{C}}		
\newtheorem{theoremA}{Theorem}
\renewcommand{\S}{\mathbb{S}}					% Sphere
\newcommand{\D}{\mathbb{D}}						% Disk

\newcommand{\dd}								% Differential d
{\mathop{}\!\mathrm{d}}						
\newcommand{\ddn}[1]							% Powers of a differential d
{\mathop{}\!\mathrm{d^{#1}}}

\newcommand{\abs}[1]							% Absolute value
{\left| #1 \right|}
\newcommand{\smallabs}[1]						% Small absolute value bars which won't scale to the argument.
{\lvert #1 \rvert}	
\newcommand{\norm}[1]							% Norm 
{\left\lVert #1 \right\rVert}	
\newcommand{\smallnorm}[1]						% Small norm bars which won't scale to the argument.
{\lVert #1 \rVert}						
\newcommand{\ip}[2]								% Inner product
{\left< #1 , #2 \right>}

\DeclareMathOperator{\diam}{diam}				
		
\DeclareMathOperator{\dist}{dist}

\newcommand{\cl}[1]{\overline{#1}}
\newcommand{\dI}{d_I}

\renewcommand{\phi}{\varphi}

\def\XXint#1#2#3{{\setbox0=\hbox{$#1{#2#3}{\int}$}
		\vcenter{\hbox{$#2#3$}}\kern-.5\wd0}}

\newtheorem{thm}{Theorem}[section]{\bf}{\it}
\newtheorem{lem}[thm]{Lemma}
\newtheorem{prop}[thm]{Proposition}
\newtheorem{cor}[thm]{Corollary}

{\bf}{\it}

{\bf}{\it}

{\bf}{\it}

{\bf}{\it}

\theoremstyle{definition}
\newtheorem{defn}[thm]{Definition}

\newtheorem{rem}[thm]{Remark}

\numberwithin{equation}{section}

\usepackage{amsthm}
\usepackage{xcolor}

\newtheorem*{customleminner}{\customlemname}

\begin{document}
\title[Homeomorphic Sobolev Extensions and the inner distance]{Homeomorphic Sobolev Extensions to Jordan Domains and the inner distance}

\author{Pekka Koskela, Chengxi Li, Xilin Zhou, Yizhe Zhu}

\address{P.\ Koskela: Department of Mathematics and Statistics, University of Jyv\"askyl\"a, P.O. Box 35, FI-40014, University of Jyv\"askyl\"a, Finland.}
\email{pekka.j.koskela@jyu.fi}
\address{C.\ Li: Department of Mathematics and Statistics, University of Jyv\"askyl\"a, P.O. Box 35, FI-40014, University of Jyv\"askyl\"a, Finland.} 
\email{cheli@jyu.fi }
\address{X.\ Zhou: Department of Mathematics and Statistics, University of Jyv\"askyl\"a, P.O. Box 35, FI-40014, University of Jyv\"askyl\"a, Finland. }
\email{zhouxilinhk2026@gmail.com}
\address{Y.\ Zhu: Institute of Mathematics, Academy of Mathematics and Systems Science, Chinese
Academy of Sciences, Beijing 100190, China.}
\email{yizhe9697@gmail.com }

\thanks{P.K. was supported by the Research Council of Finland
via Centre of Excellence in Analysis and Dynamics Research, project number 323960.
C.L. was supported by the Research Council of Finland, project numbers 364210 and 360505, and China Scholarship Council Fellowship, project number 202306200022. Y.Z.\ was supported by the Emil Aaltonen Foundation project ``\emph{Pohjois-Suomen kvasimaailman verkosto}''}

\begin{abstract}
    Given a Jordan domain $\Omega$, we characterize the situation where a given parametrization (or all parametrizations) $\varphi:\partial\D\to\partial\Omega$ has (have) a homeomorphic $W^{1,1}$-extension $f_{\varphi}\colon\overline{\D}\to\overline{\Omega}$.
\end{abstract}
\maketitle

\section{Introduction}
Let $\Omega\subset\mathbb C$ be a Jordan domain. If $\varphi\colon\partial\D \to \partial\Omega$ is a homeomorphism, then 
$\varphi$ admits a homeomorphic extension $f_{\varphi}\colon\overline{\D}\to\overline{\Omega}$. This extension need not have any regularity in general and hence it is natural to inquire for conditions under which some regularity can be gained. We say that $f_{\varphi}$ is a \emph{homeomorphic $W^{1,p}$-extension, $1\le p\le \infty$}, if $f_{\varphi}\in W^{1,p}(\D;\Omega)$ and $f_{\varphi}$ is a homeomorphic extension of $\varphi$ to $\overline{\D}$. Very recently, a rather complete solution to this problem was given in \cite{KOX25}.

\begin{theoremA}\label{thm:A}
Let $1<p<\infty$, let $\Omega$ be a Jordan domain, and let
$\varphi\colon\partial\mathbb D\to \partial\Omega$ be a homeomorphism. Then $\varphi$
admits a homeomorphic $W^{1,p}$-extension precisely when it satisfies the internal $p$-Douglas condition, i.e.,
\begin{equation*}\label{eq:internal-p-Douglas}
    \int_{\partial\mathbb D}\int_{\partial\mathbb D}
    \frac{\bigl[d_I(\varphi(x),\varphi(y))\bigr]^p}
         {|x-y|^p}\,d\mathcal{H}^1(x)\,d\mathcal{H}^1(y)
    <\infty.
\end{equation*}
\end{theoremA}
Here $\dI(x,y)$ denotes the internal distance between
$x,y\in\cl{\Omega}$, defined by
   $$ \dI(x,y):=\inf_{\gamma}\ell(\gamma),$$
where the infimum is taken over all rectifiable curves
$\gamma\subset {\Omega}\cup \{x,y\}$ joining $x$ to $y$, and
$\ell(\gamma)$ denotes the length of $\gamma$.

Theorem \ref{thm:A} leaves out the cases $p=1$ and $p=\infty$. A homeomorphism
$f:\cl{\D}\to\cl{\Omega}$ belongs to $W^{1,\infty}(\D;\Omega)$
precisely when it is Lipschitz continuous. Together with Kovalev's extension theorem
\cite[Theorem~1.1]{Kovalev2019}, the homeomorphic
$W^{1,\infty}$-extendability of $\varphi$ is equivalent to its Lipschitz
continuity. 

Our first result completes the picture by covering the case $p=1$.
\begin{thm}\label{thm:intro-criterion}
Let $\Omega$ be a Jordan domain and let
$\varphi:\partial\D\to\partial\Omega$ be a homeomorphism.
Then the following conditions are equivalent:
\begin{enumerate}
    \item[\textup{(i)}]
    $\varphi$ admits a homeomorphic $W^{1,1}$-extension.

    \item[\textup{(ii)}]
    $\varphi$ admits a homeomorphic $W^{1,1}$-extension
    $f_\varphi$ with $f_\varphi^{-1}\in W^{1,1}(\Omega;\D)$.

    \item[\textup{(iii)}]
       $\displaystyle \int_{\partial\D}d_I\bigl(\varphi(\xi),x_0\bigr)\,d\mathcal{H}^1(\xi)<\infty$
    for some (every) $x_0\in\Omega$.
\end{enumerate}
\end{thm}
The above conditions are further equivalent to the existence of a homeomorphic extension
$f_\varphi\in BV(\D;\Omega)$ with $f^{-1}_\varphi\in BV(\Omega;\D)$; see \cite[Theorem 1.3]{pratelliRadici2018} and \cite[Theorem~1.1 and its proof]{HKO2007}.
One should view the double integral in Theorem~\ref{thm:A}
as an intrinsic trace norm for homeomorphic $W^{1,p}$-mappings.
Our Theorem~\ref{thm:intro-criterion} then gives an intrinsic
homeomorphic analogue for the fact that the trace space of the
usual $W^{1,1}(\D;\Omega)$ is $L^1(\partial\D)$.

By \cite{Verchota2007}, every homeomorphism $\varphi\colon\partial\D\to \partial\D$
has a homeomorphic $W^{1,p}$-extension
$f_\varphi\colon\cl{\D}\to\cl{\D}$ when $1\le p<2$, but not
necessarily when $p\ge2$. This, together with the wish to
understand when the natural energy minimization problem makes
sense for arbitrary homeomorphic boundary data, has motivated
studies on conditions on $\Omega$ towards such extendability
\cite{BJKXZ25,KKO20,KO21,KO23Bi,KO23JS,KOX25,KWX18}. No actual characterizations are known.

Our second result gives a full characterization of those Jordan
domains for which all boundary parametrizations have
homeomorphic $W^{1,1}$-extensions. We define the \emph{inner diameter} of $\Omega$ by
$$\diam_I(\Omega)
    :=\sup\{\dI(x,y):x,y\in\Omega\}.$$

\begin{thm}\label{thm:intro-all-boundary-data}
Let $\Omega\subset\mathbb C$ be a Jordan domain.
Then every homeomorphism $\varphi:\partial\D\to\partial\Omega$
admits a homeomorphic $W^{1,1}$-extension if and only if
$\diam_I(\Omega)<\infty$.
\end{thm}

Let us single out a couple of consequences of
Theorem~\ref{thm:intro-all-boundary-data} and its proof.

\begin{cor}\label{cor:uniform-W11-extension}
Let $\Omega\subset\mathbb C$ be a Jordan domain.
The following conditions are equivalent:
\begin{enumerate}
    \item[\textup{(i)}]
    Every boundary homeomorphism
    $\varphi:\partial\mathbb D\to\partial\Omega$
    admits a homeomorphic $W^{1,1}$-extension.

    \item[\textup{(ii)}]
    There exists a constant $C_\Omega<\infty$,
    depending only on $\Omega$, such that every boundary
    homeomorphism $\varphi:\partial\mathbb D\to\partial\Omega$
    admits a homeomorphic bi-Sobolev extension
    $f_{\varphi}:\overline{\mathbb D}\to\overline\Omega$ satisfying $$ \|f_{\varphi}\|_{W^{1,1}(\mathbb D;\Omega)}+\|f^{-1}_{\varphi}\|_{W^{1,1}(\Omega;\D)}\le C_\Omega.$$
\end{enumerate}
\end{cor}

\begin{cor}\label{cor:intro-sharpness}
There exist a Jordan domain $\Omega$ and a constant $M<\infty$,
such that every homeomorphism $\varphi:\partial\D\to\partial\Omega$
admits a homeomorphic $W^{1,1}$-extension $f_\varphi$ with
\[
    \|f_\varphi\|_{W^{1,1}(\D;\Omega)}\le M,
\]
but there exists a homeomorphism $\varphi\colon\partial\D\to\partial\Omega$
that admits no homeomorphic $W^{1,p}$-extension for any $p>1$.
\end{cor}

In the case of usual Sobolev extensions, the existence of an extension
for every function implies a uniform norm bound in terms
of the original norm by the open mapping theorem; see
\cite[Section~3, p.~1224]{HKT2008}. In our setting, this argument does not apply.
Corollary~\ref{cor:intro-sharpness} extends a recent result
from \cite[Theorem~1.3 and Example~3.1]{KoskelaXu2026} to the missing case $p=1$.

We close this introduction with the following corollary.

\begin{cor}\label{thm: equivalent-extension-properties}
Let $\Omega\subset\mathbb C$ be a Jordan domain.
The following conditions are equivalent, where each of
\textup{(i)}--\textup{(v)} is understood to hold for every orientation-preserving boundary homeomorphism
$\varphi:\partial\D\to\partial\Omega$:
\begin{enumerate}
    \item[\textup{(i)}]
    $\varphi$
    admits a homeomorphic extension of finite distortion.

    \item[\textup{(ii)}]
    $\varphi$
    admits a homeomorphic bi-Sobolev extension.

    \item[\textup{(iii)}]
    $\varphi$
    admits a homeomorphic $W^{1,1}$-extension.

    \item[\textup{(iv)}]
    $\varphi$
    admits a homeomorphic $BV$-extension.

     \item[\textup{(v)}]
    $\varphi$
    admits a homeomorphic bi-$BV$-extension.

    \item[\textup{(vi)}]
    $\operatorname{diam}_I(\Omega)<\infty$.
\end{enumerate}
\end{cor}

Here, a homeomorphic extension
$f_{\varphi}\colon\cl{\D}\to\cl{\Omega}$ is said to be of
\emph{finite distortion} if
$f_{\varphi}\in W^{1,1}(\D;\mathbb C)$,
$J_{f_{\varphi}}\in L^1(\D)$, and
\[
    |Df_{\varphi}(z)|^2\le K(z)J_{f_{\varphi}(z)}
    \qquad\text{for a.e. }z\in\D,
\]
for some measurable function $K:\D\to[1,\infty)$.
Here $J_{f_{\varphi}}$ denotes the Jacobian determinant of $f$,
and $|Df_{\varphi}|$ denotes the operator norm of its derivative.
We call $f_{\varphi}$ a \emph{bi-Sobolev extension} if
\[
    f_{\varphi}\in W^{1,1}(\D;\Omega)
    \quad\text{and}\quad
    f_{\varphi}^{-1}\in W^{1,1}(\Omega;\D).
\]
A \emph{bi-$BV$ extension} is defined analogously,
with $W^{1,1}$ replaced by $BV$.

In fact, some implications in Corollary \ref{thm: equivalent-extension-properties} are trivial and some follow from earlier results.
Indeed, the implications
$\textup{(i)}\Rightarrow\textup{(iii)}$,
$\textup{(ii)}\Rightarrow\textup{(iii)}$, and
$\textup{(v)}\Rightarrow\textup{(iv)}$
are immediate from the definitions.
Furthermore, the inclusion $W^{1,1}\subset BV$ yields
$\textup{(iii)}\Rightarrow\textup{(iv)}$ and
$\textup{(ii)}\Rightarrow\textup{(v)}$.

Next, the implication \textup{(i)}$\Rightarrow$\textup{(ii)}
follows from ~\cite[Theorem~1.2]{HenclKoskela2006}.
Indeed, if $f\in W^{1,1}(\D;\Omega)$ is a homeomorphism
of finite distortion, then
$f^{-1}\in W^{1,1}_{\mathrm{loc}}(\Omega;\D)$ and
\[
    \int_\Omega |D(f^{-1})(y)|\,dA(y)
    =\int_\D |Df(x)|\,dA(x)<\infty.
\]
Since $\Omega$ is bounded and $f^{-1}(\Omega)=\D$,
we conclude that $f^{-1}\in W^{1,1}(\Omega;\D)$.
The reverse implication \textup{(ii)}$\Rightarrow$\textup{(i)} follows from
\cite[Theorem~1.6]{HenclKoskela2014}.

The implication \textup{(iv)}$\Rightarrow$\textup{(v)}
follows from \cite[Theorem~1.1 and its proof]{HKO2007}. Furthermore, \cite[Theorem~1.3]{pratelliRadici2018}
yields \textup{(iv)}$\Rightarrow$\textup{(iii)}:
the existence of a homeomorphic $BV$-extension guarantees
the existence of a homeomorphic $W^{1,1}$-extension
with the same boundary values.

For the convenience of the reader, we provide essentially
self-contained proofs for our results in what follows.

In Section \ref{preliminary}, we recall some basic notations and collect some lemmas which will be useful in our main proofs. In Section \ref{circleextension} we construct a homeomorphic $W^{1,1}$-extension for every circle homeomorphism; this is the most technical contributions in our paper. We prove Theorem \ref{thm:intro-all-boundary-data} and Corollaries \ref{cor:uniform-W11-extension}-\ref{thm: equivalent-extension-properties} in Section \ref{thm1.2} and prove Theorem \ref{thm:intro-criterion} in Section \ref{thm1.1}.

\section{Preliminaries}\label{preliminary}

\subsection{Notation}
Throughout this paper, a {\em Jordan domain} means a bounded domain in $\mathbb C$ whose boundary is a Jordan curve.

\subsection{Hyperbolic metric}
In this subsection, we recall the definition of the hyperbolic metrics and collect some basic properties needed below.
For a detailed discussion, we
refer to the monograph of Astala, Iwaniec, and Martin
\cite[Chapter~2]{AIM09}.

\begin{defn}\label{def:hyperbolic-metric}
For $z,w\in\D$, \emph{the hyperbolic metric} between $z$ and $w$ in $\D$ is defined by
\[
    d_{\D}(z,w)
    :=\inf_{\gamma}\int_{\gamma}
    \frac{2}{1-|y|^2}\,ds(y),
\]
where the infimum is taken over all rectifiable curves
$\gamma\subset\D$ joining $z$ to $w$.

Let $\Omega\subsetneq\mathbb C$ be a simply connected domain,
and let $g\colon\D\to\Omega$ be a conformal homeomorphism.
The hyperbolic metric in $\Omega$ is defined by
\[
    d_{\Omega}(x,y)
    :=d_{\D}\bigl(g^{-1}(x),g^{-1}(y)\bigr),
    \qquad x,y\in\Omega.
\]
This definition is independent of the choice of $g$, since every
conformal automorphisms of $\D$ is an isometry with respect to $d_{\D}$.
\end{defn}

Since conformal homeomorphisms preserve hyperbolic distances and radial segments
are hyperbolic geodesics in $\D$ (see \cite[p.~2]{Ahl73}), the following lemma is immediate.
\begin{lem}\label{lem:conformal-invariance}
Let $g\colon\D\to\Omega$ be a conformal homeomorphism.
Then, for every $0<R<1$ and $t\in\mathbb R$,
the curve $g([0,Re^{it}])$ is a hyperbolic geodesic segment
in $\Omega$ joining $g(0)$ to $g(Re^{it})$.
\end{lem}

We will also use the following form of the Koebe distortion theorem;
see \cite[Theorem~2.10.8]{AIM09}.

\begin{lem}\label{lem: koebe ditortion theorem}
Let $z,w\in\D$ and let $\varphi\colon\mathbb D\to\Omega\subset\mathbb C$
be a conformal map. Then
$$e^{-3d_{\D}(z,w)}\le\frac{|\varphi'(z)|}{|\varphi'(w)|}\le e^{3d_{\D}(z,w)}.$$
\end{lem}

% We next introduce the inner distance.
% \begin{defn}\label{defofinner}
% Let $\Omega\subset\mathbb{R}^2$ be a Jordan domain. For
% $x,y\in\overline{\Omega}$, we define \emph{the inner distance} by
% \[
%     d_I(x,y):=\inf_{\gamma}\ell(\gamma),
% \]
% where the infimum is taken over all continuous curves
% $\gamma:[0,1]\to\mathbb{R}^2$ joining $x$ to $y$ with $\gamma((0,1))\subset\Omega$.
% Here $\ell(\gamma)\in[0,+\infty]$ denotes the Euclidean length
% of $\gamma$. In particular, $d_I(x,y):=+\infty$ if no such
% curve is rectifiable.
% \end{defn}

Finally, we recall the Gehring--Hayman inequality, which compares
the Euclidean length of a hyperbolic geodesic with the inner
distance between its endpoints. For a proof, see \cite{GH62}
and \cite[Theorem~4.20]{Pom92}.

\begin{lem}
\label{lem: gehring-hayman}
Let $\Omega\subsetneq\mathbb C$ be a simply connected
domain. Given  $x,y\in\Omega$, let $\gamma$ be the hyperbolic geodesic joining $x$ to $y$ and let $\sigma$ be any other curve in $\Omega$ connecting them. Then 
\[
    \ell(\gamma)\le C\,\ell(\sigma)
\]
for constant $C$ independent of $x,y,\Omega$.
\end{lem}

\subsection{Preliminary results} In this subsection, we collect some necessary results which will be useful in our main proofs. We first give the polar-coordinate version of the Sobolev and BV
slicing theorems. It follows by a change of variables from
\cite[Section~3.11]{ambrosioFuscoPallara2000}.
\begin{thm}
\label{thm:sobolev-slicing}
Let $c\in\mathbb C$, $0<a<b<\infty$, and
$\alpha<\beta$ with $\beta-\alpha<2\pi$. Set
\[
    S:=\{c+re^{it}:a<r<b,\ \alpha<t<\beta\},
\]
and let $u\colon\overline S\to\mathbb R^2$ be continuous.
We define the radial and angular slices by
$$ u_t(r):=u(c+re^{it}),\  a\le r\le b\quad \text{and}\quad v_r(t):=u(c+re^{it}), \ \alpha\le t\le\beta,$$
respectively.

\begin{enumerate}
    \item[\textup{(i)}]
    Let $u\in W^{1,p}(S;\mathbb R^2)$ for some
    $1\le p<\infty$. Then for a.e.
    $t\in(\alpha,\beta)$, the slice $u_t$ belongs to
    $W^{1,p}((a,b);\mathbb R^2)$, is absolutely continuous
    on $[a,b]$, and satisfies
    \[
        u_t'(r)=Du(c+re^{it})e^{it}
        \quad\text{for a.e. }r\in(a,b).
    \]
    Similarly, for a.e. $r\in(a,b)$, the slice
    $v_r$ belongs to $W^{1,p}((\alpha,\beta);\mathbb R^2)$,
    is absolutely continuous on $[\alpha,\beta]$, and satisfies
    \[
        v_r'(t)=rDu(c+re^{it})(ie^{it})
        \quad\text{for a.e. }t\in(\alpha,\beta).
    \]

    \item[\textup{(ii)}]
    Let $u\in BV(S;\mathbb R^2)$. Then for a.e.
    $t\in(\alpha,\beta)$, the slice $u_t$ belongs to
    $BV((a,b);\mathbb R^2)$ and
    \[
        \ell(u_t)=|Du_t|((a,b)).
    \]
    Moreover,
    \[
        \int_\alpha^\beta \ell(u_t)\,dt
        \le \frac{1}{a}|Du|(S).
    \]
    Here $Du_t$ denotes its distributional
derivative with respect to $r$, and $|Du_t|$ denotes the associated total variation measure.
\end{enumerate}
\end{thm}

The following technical lemma shows that a continuous, strictly
increasing function $\eta$ satisfying
$\eta(t+2\pi)=\eta(t)+2\pi$ can be uniformly approximated by
smooth functions with the same properties. We include a proof for completeness. 

\begin{lem}\label{lem: approxmiate by convex combination}
Let $\eta\colon\mathbb R\to\mathbb R$ be a continuous,
strictly increasing function satisfying $\eta(t+2\pi)=\eta(t)+2\pi$.
Then there exists a sequence of smooth functions
$\eta_n\colon\mathbb R\to\mathbb R$ such that
\begin{equation}\label{ineq: conditions of eta_n}
     \eta_n(t+2\pi)=\eta_n(t)+2\pi,
    \quad
    \|\eta_n-\eta\|_{\infty}\le2^{-n},
    \quad\text{and}\quad
     \eta_n'>\frac{2^{-n}}{4\pi}.
\end{equation}
\end{lem}

\begin{proof}
Choose a nonnegative function
$\rho\in C_0^\infty((-1,1))$ with
$\int_{\mathbb R}\rho=1$, and set
\[
    \rho_\varepsilon(t):=\varepsilon^{-1}\rho(t/\varepsilon)
\]
for every $\varepsilon>0$.
Since the function $t\mapsto\eta(t)-t$ is continuous and periodic,
$\eta$ is uniformly continuous on $\mathbb R$.
Thus, for every $n\in\N$, we may choose
$\varepsilon_n>0$ such that
\[
    \sup_{\substack{t\in\mathbb R,\, |s|<\varepsilon_n}}
    |\eta(t-s)-\eta(t)|
    \le 2^{-(n+1)}.
\]
Define
\[
    u_n(t):=(\eta*\rho_{\varepsilon_n})(t)
    =\int_{\mathbb R}\eta(t-s)\rho_{\varepsilon_n}(s)\,ds.
\]
Then $u_n$ is smooth with
$u_n'>0$.
Moreover,
\[
    u_n(t+2\pi)=u_n(t)+2\pi
    \quad\text{and}\quad
    \|u_n-\eta\|_{L^\infty(\mathbb R)}
    \le2^{-(n+1)}.
\]

We now define the sequence $\eta_n$ by
\[
    \eta_n(t):=(1-\delta_n)u_n(t)+\delta_n a(t),
\]
where $\delta_n=2^{-n}/4\pi$ and $a(t):=t+\eta(0)$. Clearly, $\eta-a$ is $2\pi$-periodic, and 
\[
    \|\eta-a\|_{L^\infty(\mathbb R)}=\|\eta-a\|_{L^\infty([0,2\pi])}\le2\pi.
\]
By construction, each $\eta_n$ is smooth and satisfies
\[
    \eta_n(t+2\pi)=\eta_n(t)+2\pi.    
\]
Moreover, since $u_n'>0$ and $a'=1$, we have
\[
    \eta_n'(t)
    =(1-\delta_n)u_n'(t)+\delta_n
    \ge\delta_n.
\]
Finally, the choice of $\delta_n$ gives
\[
    \|\eta_n-\eta\|_{L^\infty(\mathbb R)}
    \le
    (1-\delta_n)\|u_n-\eta\|_{L^\infty(\mathbb R)}
    +\delta_n\|a-\eta\|_{L^\infty(\mathbb R)}\\
    \le 2^{-n}.
\]
This completes the proof.
\end{proof}

\section{Homeomorphic extensions of circle maps}\label{circleextension}

In this section, we construct a suitable homeomorphic $W^{1,1}$-extension
for a given homeomorphism $\psi\colon\partial\D\to\partial\D$.
This construction is one of the main technical contributions
of the paper. 

Our extension is locally
Lipschitz in $\D$ with estimates controlling both its
radial and angular behavior; see Proposition \ref{prop:circle-extension}. Moreover, its $W^{1,1}$-norm is bounded by an absolute constant;
see Proposition \ref{prop:circle-extension-W11}.

\begin{prop}\label{prop:circle-extension}
Let $\psi\colon\partial\D\to\partial\D$ be an orientation-preserving
homeomorphism. Then $\psi$ admits an orientation-preserving homeomorphic
extension $\Psi\colon\overline{\D}\to\overline{\D}$ which is locally
bi-Lipschitz in $\D$ with
\[
    |\Psi(z)|=|z| \quad \text{in }\overline{\D},
    \quad
    |\partial_r\Psi(re^{it})|\le 10 \qquad \text{a.e. in }\D,
\]
and
\[
    d_{\D}\bigl(\Psi(re^{it}),r\psi(e^{it})\bigr)\le 4
    \qquad \text{for every }0<r<1\text{ and }t\in\mathbb R.
\]
Moreover, for some $\alpha\in\mathbb R$,
\[
    \Psi(z)=e^{i\alpha}z
    \qquad \text{for }|z|\le \tfrac12.
\]
\end{prop}

\begin{proof}
Let $\eta\colon\mathbb R\to\mathbb R$ be a continuous, strictly
increasing lift of $\psi$, i.e.,
\[
    \psi(e^{it})=e^{i\eta(t)},
    \qquad
    \eta(t+2\pi)=\eta(t)+2\pi.
\]
By Lemma~\ref{lem: approxmiate by convex combination}, there exists
a sequence of smooth functions $\eta_n\colon\mathbb R\to\mathbb R$,
$n\ge 2$, such that
\begin{equation}\label{propofetan}
    \eta_n(t+2\pi)=\eta_n(t)+2\pi,
    \quad
    \eta_n'(t)\ge \frac{2^{-n}}{4\pi},
    \quad \text{and} \quad
    \|\eta_n-\eta\|_\infty\le 2^{-n}.
\end{equation}

We choose
\[
    \alpha:=\frac12\left(
        \max_{t\in[0,2\pi]}\bigl(\eta_2(t)-t\bigr)
        +\min_{t\in[0,2\pi]}\bigl(\eta_2(t)-t\bigr)
    \right)
\]
and set 
\[
    R_\alpha(t):=t+\alpha.
\]
Since $\eta_2$ is increasing and satisfies
$\eta_2(t+2\pi)=\eta_2(t)+2\pi$, the function $t\mapsto\eta_2(t)-t$
is $2\pi$-periodic with oscillation at most $2\pi$. Hence
\[
    \|\eta_2-R_\alpha\|_\infty\le\pi.
\]

Let $r_n:=1-2^{-n}$ for $n\ge 2$. We define
\begin{equation}\label{defofArt}
    A(r,t):=
    \begin{cases}
        R_\alpha(t),
            &0\le r\le \tfrac12,\\[1mm]
        (3-4r)R_\alpha(t)+(4r-2)\eta_2(t),
            &\tfrac12\le r\le \tfrac34,\\[1mm]
        (1-\lambda_n(r))\eta_n(t)+\lambda_n(r)\eta_{n+1}(t),
            &r_n\le r\le r_{n+1},
    \end{cases}
\end{equation}
where $\lambda_n(r):=(r-r_n)/(r_{n+1}-r_n)$.
The formulas agree at their common endpoints. For
$r_n\le r\le r_{n+1}$, we have
\begin{equation}\label{eq:angular-approximation}
    |A(r,t)-\eta(t)|
    \le (1-\lambda_n(r))\|\eta-\eta_n\|_\infty
        +\lambda_n(r)\|\eta-\eta_{n+1}\|_\infty
    \le 1-r.
\end{equation}
The function $A$ is piecewise affine in $r$. In particular, the
partial derivative $A_r:=\partial_rA$ satisfies
\[
    |A_r(r,t)|
    \le \frac{\|\eta_{n+1}-\eta_n\|_\infty}{r_{n+1}-r_n}
    \le \frac{2^{-n-1}+2^{-n}}{2^{-n-1}}=3
    \qquad \text{for }r_n<r<r_{n+1}.
\]
For $1/2<r<3/4$, we have
\[
    |A_r(r,t)|=4|\eta_2(t)-R_\alpha(t)|\le 4\pi,
\]
and $A_r=0$ for $0<r<1/2$. Consequently,
\begin{equation}\label{ineq: upper bound of A_r}
    |A_r(r,t)|\le 4\pi
    \qquad \text{for a.e. }(r,t)\in(0,1)\times\mathbb R.
\end{equation}

Moreover, by \eqref{propofetan} and \eqref{defofArt}, we have
\begin{equation}\label{periodicityofA}
    A(r,t+2\pi)=A(r,t)+2\pi,
\end{equation}
and the partial derivative $A_t:=\partial_t A$ satisfies
\[
    A_t(r,t)\ge \frac{2^{-n}}{8\pi}
    \qquad \text{for }r_n\le r\le r_{n+1}.
\]
For $1/2\le r\le 3/4$, the derivative $A_t$ is a convex
combination of $1$ and $\eta_2'(t)$, so $A_t\ge 1/(16\pi)$.
Since $A_t=1$ for $0\le r\le 1/2$, we obtain
\begin{equation}\label{ineq: A_t}
    A_t(r,t)\ge \frac{1-r}{8\pi}>0
    \qquad \text{for }0\le r<1\text{ and }t\in\mathbb R.
\end{equation}

Now we define the extension of $\psi$ by
\begin{equation}\label{defofPsi}
    \Psi(re^{it}):=re^{iA(r,t)}
    \quad 0\le r<1,
    \qquad
    \Psi|_{\partial\D}:=\psi.
\end{equation}
By \eqref{periodicityofA}, this definition is well defined.
We first show that $\Psi$ is a sense-preserving homeomorphism of
$\overline{\D}$. It is continuous in $\D$ by construction and at
the boundary by \eqref{eq:angular-approximation}. Since $\Psi$
fixes the origin and maps each circle $r\partial\D$
homeomorphically onto itself, it is bijective and hence a
homeomorphism. Moreover, $\Psi$ is sense-preserving.

We next show that $\Psi$ is locally bi-Lipschitz in $\D$.
Since $\Psi$ is a rotation on $B(0,1/2)$, it suffices to
consider points away from the origin. There, $\Psi$ is locally
Lipschitz by the local Lipschitz continuity of $A$. For each
$0<r<1$, let $A^{-1}(r,\cdot)$ denote the inverse of
$A(r,\cdot)$ by $A(r,A^{-1}(r,\theta))=\theta$. Since $\Psi$ preserves the modulus, we have
\[
    \Psi^{-1}(re^{i\theta})=re^{iA^{-1}(r,\theta)}.
\]
On each annulus  in \eqref{defofArt}, the positivity of
$A_t$ gives the smoothness of $A^{-1}$ by the implicit function
theorem. Applying the chain rule, we obtain
\[
    \partial_\theta A^{-1}(r,\theta)
    =\frac{1}{A_t(r,A^{-1}(r,\theta))},
    \qquad
    \partial_r A^{-1}(r,\theta)
    =-\frac{A_r(r,A^{-1}(r,\theta))}
            {A_t(r,A^{-1}(r,\theta))}.
\]
for a.e. $(r,t)\in(0,1)\times\R$.

Fix $0<R<1$ and set $m_R:=(1-R)/(8\pi)$. By
\eqref{ineq: upper bound of A_r} and \eqref{ineq: A_t},
we have $|A_r|\le 4\pi$ and $A_t\ge m_R$ for a.e. $0<r<R$. Hence
\[
\begin{aligned}
    |D\Psi^{-1}(re^{i\theta})|
    &\le |\partial_r\Psi^{-1}(re^{i\theta})|
        +\frac1r|\partial_\theta\Psi^{-1}(re^{i\theta})|\\
    &\le 1+\frac{4\pi R+1}{m_R}.
\end{aligned}
\]
Together with the continuity of $\Psi^{-1}$, this estimate
implies that $\Psi^{-1}$ is Lipschitz on $B(0,R)$. 
Since $R$ is arbitrary, $\Psi$ is locally bi-Lipschitz in $\D$.

By the piecewise estimates for $|A_r|$ established above, we obtain
\[
    r|A_r(r,t)|\le
    \begin{cases}
        0,    &0\le r<1/2,\\
        3\pi, &1/2\le r<3/4,\\
        3,    &3/4\le r<1
    \end{cases}
\]
for a.e. $(r,t)\in[0,1]\times\R$.
Consequently,
\begin{equation}\label{ineq: upper bound of partial_r Psi}
    |\partial_r\Psi(re^{it})|
    =\sqrt{1+r^2|A_r(r,t)|^2}
    \le \sqrt{1+9\pi^2}<10
    \qquad \text{a.e. in }\D.
\end{equation}

Finally, we estimate the hyperbolic distance between $\Psi(re^{it})$ and $r\psi(e^{it})$. For $3/4\le r<1$, let
\[
    \gamma_r(s):=re^{i[(1-s)A(r,t)+s\eta(t)]}
    \qquad \text{for }0\le s\le 1,
\]
i.e., $\gamma_r$ is the circular arc of radius $r$ that joins
$\Psi(re^{it})$ to $r\psi(e^{it})$. By
Definition~\ref{def:hyperbolic-metric} and
\eqref{eq:angular-approximation}, we conclude that
\[
    d_{\D}\bigl(\Psi(re^{it}),r\psi(e^{it})\bigr)
    \le \ell_{\D}(\gamma_r)
    \le \frac{2r}{1-r^2}|A(r,t)-\eta(t)|
    \le 1.
\]
For $0<r\le 3/4$, let $\gamma$ be the path that follows the
radial segment from $\Psi(re^{it})$ to $0$ and then the radial
segment from $0$ to $r\psi(e^{it})$. Since both endpoints of
$\gamma$ have modulus $r$, we have
\[
\begin{aligned}
    d_{\D}\bigl(\Psi(re^{it}),r\psi(e^{it})\bigr)
    &\le \ell_{\D}(\gamma)
     = 2\int_0^r \frac{2\,ds}{1-s^2}\\
    &= 2\log\frac{1+r}{1-r}
     \le 2\log 7 < 4.
\end{aligned}
\]
This completes the proof.
\end{proof}

We now estimate the $W^{1,1}$-norm of the homeomorphic extension constructed above.

\begin{prop}\label{prop:circle-extension-W11}
Let $\psi\colon\partial\D\to\partial\D$ be an orientation-preserving
homeomorphism. The homeomorphic extension
$\Psi:\overline{\D}\to\overline{\D}$ constructed in
Proposition~\ref{prop:circle-extension} belongs to the Sobolev class
$W^{1,1}(\D;\D)$.
\end{prop}

\begin{proof}
Since $|\Psi|\in L^1(\D)$, it remains to prove that  $|D\Psi|\in L^1(\D)$. By the construction of $\Psi$ given in Proposition \ref{prop:circle-extension}, $\Psi$ is locally Lipschitz in $\D$. The polar coordinate yields that 
\begin{equation}\label{ineq: polar coordinates}
    \int_{\D}|D\Psi(z)|\,dA(z)
    \le \int_0^1\int_0^{2\pi}
    \bigl(r|\partial_r\Psi(re^{it})|
    +|\partial_t\Psi(re^{it})|\bigr)
    \,dt\,dr.
\end{equation}

The finiteness of the first term on the right of \eqref{ineq: polar coordinates}
follows directly from \eqref{ineq: upper bound of partial_r Psi}. 
For the second term, the map $t\mapsto\Psi(re^{it})$
is Lipschitz on $[0,2\pi]$, maps $[0,2\pi)$ bijectively
onto $\S(0,r)$, and has equal endpoint values.
Thus, by the arc length formula, we have
    $$\int_0^{2\pi}|\partial_t\Psi(re^{it})|\,dt
    =\ell(\S(0,r))=2\pi r.$$
and we conclude that
$$ \|\Psi\|_{W^{1,1}(\D;\D)}=\int_{\D}|\Psi(z)|\,dA(z)+\int_{\D}|D\Psi(z)|\,dA(z)\le 6\pi. $$ 
 This completes the proof.
\end{proof}

\section{The case of all parametrizations}\label{thm1.2}
In this section, we prove one of our main results,
Theorem~\ref{thm:intro-all-boundary-data}, together with
Corollaries~\ref{cor:uniform-W11-extension}--%
\ref{thm: equivalent-extension-properties}. 

\subsection{The proof of Theorem \ref{thm:intro-all-boundary-data}}
We prove Theorem~\ref{thm:intro-all-boundary-data} in two parts. For sufficiency, we show that if $\diam_I(\Omega)<\infty$, then every boundary homeomorphism
$\varphi\colon\partial\D\to\partial\Omega$ has a homeomorphic $W^{1,1}$-extension; see Proposition~\ref{thm: W11-extension}. For necessity, we obtain a stronger result: if
$\diam_I(\Omega)=\infty$, then there exists a boundary homeomorphism
that does not even have a homeomorphic $BV$-extension; see
Proposition~\ref{main-3}. 

We begin with the proof of sufficiency.
\begin{prop}\label{thm: W11-extension}
    Let $\Omega\subset\C$ be a Jordan domain with $M_\Omega:=\diam_I(\Omega)<\infty$. Then every boundary homeomorphism $\varphi\colon\partial\D\to\partial\Omega$ admits a homeomorphic extension $f_\varphi:\overline{\D}\to\overline{\Omega}$ in the Sobolev class $W^{1,1}(\D;\Omega)$. Moreover, we have 
    $$\int_{\D}|Df_\varphi(z)|\,dA(z)\le C M_\Omega$$
  for an absolute constant $C>0$.
\end{prop}
\begin{proof}
It suffices to consider orientation-preserving boundary maps, the other case follows by precomposing both the boundary map
and its extension with the reflection $z\mapsto\overline z$.

Let $g:\D\to\Omega\subsetneq \mathbb C$ be conformal. By Carath\'eodory's theorem \cite[page 18]{Pom92},
$g$ extends to a homeomorphism
$g:\overline{\D}\to\overline{\Omega}$.
For every $t\in\R$ and $0<R<1$,
Lemma~\ref{lem:conformal-invariance} shows that
$g([0,Re^{it}])$ is a hyperbolic geodesic in $\Omega$.
Hence, Lemma \ref{lem: gehring-hayman} gives
$$    \int_0^R |g'(re^{it})|\,dr
    \le C\,d_I\bigl(g(0),g(Re^{it})\bigr)
    \le C M_\Omega.$$
Letting $R\to 1$, we obtain
\begin{equation}\label{eq:uniform-radial-length'}
    \sup_{t\in\R}
    \int_0^1 |g'(re^{it})|\,dr
    \le C M_\Omega.
\end{equation}

Set $\psi:=(g|_{\partial\D})^{-1}\circ\varphi:\partial\D \to \partial\D$, and let
$\Psi:\overline{\D}\to\overline{\D}$ be its extension given by
Proposition~\ref{prop:circle-extension}.
Then $f_\varphi:=g\circ\Psi$ is a homeomorphism from
$\overline{\D}$ onto $\overline{\Omega}$, and
locally Lipschitz in $\D$ with $f_\varphi|_{\partial\D}=\varphi$.

It remains to prove $f_\varphi\in W^{1,1}(\D;\Omega)$. Since $\operatorname{diam}(\Omega)\le \operatorname{diam}_I(\Omega)$, $\Omega$ is bounded, and we only need to show that $|Df_\varphi|\in L^1(\D)$. Note that 
\begin{multline}\label{ineq: inter DF'}
    \int_{\D} |Df_\varphi(z)|\,dA(z)=\int_0^1\int_0^{2\pi}|Df_\varphi(re^{it})|r\,dt\,dr\\
    \le \int_0^1\int_0^{2\pi}|\partial_rf_\varphi(re^{it})|r\,dt\,dr+ \int_0^1\int_0^{2\pi}|\partial_tf_\varphi(re^{it})|\,dt\,dr:=I_1+I_2.
\end{multline}

We first estimate $I_1$ in \eqref{ineq: inter DF'}.
By Proposition~\ref{prop:circle-extension}, we have
    $$d_{\D}\bigl(\Psi(re^{it}),r\psi(e^{it})\bigr)\le4.$$
The Koebe distortion theorem (Lemma \ref{lem: koebe ditortion theorem}) therefore gives
    $$|g'(\Psi(re^{it}))|\le e^{3d_{\D}\bigl(\Psi(re^{it}),r\psi(e^{it})\bigr)}|g'(r\psi(e^{it}))|
    \le e^{12}|g'(r\psi(e^{it}))|.$$
Combining this with $|\partial_r\Psi|\le10$ and the chain rule,
we obtain
    $$|\partial_r f_\varphi(re^{it})|
    \le10e^{12}|g'(r\psi(e^{it}))|
    \le C|g'(r\psi(e^{it}))| \qquad\text{a.e. in }\D.$$
Thus, by \eqref{eq:uniform-radial-length'}, Tonelli's theorem, and $r\le1$, we deduce that
\begin{equation}\label{eq:radial-energy}
   I_1
    \le C\int_0^{2\pi}\int_0^1
    |g'(r\psi(e^{it}))|\,dr\,dt\le C M_\Omega.
\end{equation}

To estimate $I_2$, fix $0<r_0<1$.
Since $\Psi$ is locally Lipschitz and maps $\S(0,r_0)$ homeomorphically onto itself
and $g$ is conformal, the maps $t\mapsto f_\varphi(r_0e^{it})$ and $ t\mapsto g(r_0e^{it})$
are Lipschitz parametrizations of $g(\S(0,r_0))$
on $[0,2\pi]$, each bijective on $[0,2\pi)$
and with equal endpoint values.
The arc length formula therefore gives
 $$\int_0^{2\pi}|\partial_tf_\varphi(r_0e^{it})|\,dt
    =\ell\bigl(g(\S(0,r_0))\bigr)
    =r_0\int_0^{2\pi}|g'(r_0e^{it})|\,dt.$$
Integrating with respect to $r_0$, combining with
\eqref{eq:uniform-radial-length'} and using Tonelli's theorem, it follows that
\begin{equation}\label{eq:angular-energy}
   I_2=\int_0^{2\pi}\int_0^1
    r|g'(re^{it})|\,dr\,dt\le C M_\Omega.
\end{equation}

Finally, combining the above estimate with \eqref{ineq: inter DF'} and \eqref{eq:radial-energy}, we conclude that
$$    \int_{\D}|Df_\varphi(z)|\,dA(z)
    \le I_1+I_2
    \le C M_\Omega<\infty.$$
This completes the proof.
\end{proof}

\begin{rem}\label{rem: bi-Sobolev}
In fact, the extension
$f_\varphi$ constructed in Proposition~\ref{thm: W11-extension} also satisfies
\begin{equation}\label{ineq: bi-sobolev energy estimate}
    \int_\Omega |Df_\varphi^{-1}(w)|\,dA(w)=\int_{\D}|Df_\varphi(z)|\,dA(z)\le CM_\Omega.
\end{equation}

Indeed, Proposition~\ref{prop:circle-extension} shows that
$\Psi$ is locally bi-Lipschitz in $\D$.
Since $g$ is conformal, $f_\varphi=g\circ\Psi$ is also
locally bi-Lipschitz.
The change of variables formula and the chain rule therefore yield
$$\begin{aligned}
   \begin{aligned}
    \int_\Omega |Df_\varphi^{-1}(w)|\,dA(w)
    &=\int_{\D}|Df_\varphi^{-1}(f_\varphi(z))||J_{f_\varphi}(z)|\,dA(z)\\
    &=\int_{\D}|(Df_\varphi(z))^{-1}||J_{f_\varphi}(z)|\,dA(z)\\
    &=\int_{\D}|Df_\varphi(z)|\,dA(z)\le CM_\Omega.
\end{aligned}
\end{aligned}$$
Since $f_\varphi^{-1}$ is bounded and $\Omega$ has finite area,
it follows that
$f_\varphi^{-1}\in W^{1,1}(\Omega;\D)$.

\end{rem}

We then turn to the proof of necessity.

\begin{prop}\label{main-3}
    Let $\Omega$ be a Jordan domain with $\diam_I(\Omega) = \infty$.
    Then there is a boundary homeomorphism  $\varphi\colon \partial\D \to \partial\Omega$ so that $\varphi$ does not have a homeomorphic BV-extension.
\end{prop}
\begin{proof}
    Given $\omega\in\partial\D$ and $0<\delta<\pi$, let
$$I(\omega,\delta):=\{\omega e^{i\alpha}:|\alpha|\le\delta\}\subset \partial\D$$
denote the closed arc of $\partial\D$ centered at $\omega$
with angular length $2\delta$.
    By the Riemann mapping theorem, there is a conformal homeomorphism $g\colon\D\to\Omega$. Since $\Omega$ is a Jordan domain, Carath\'eodory's theorem ensures
that $g$ extends to a homeomorphism
$g\colon\overline{\D}\to\overline{\Omega}$.
    
    We first show that there exists $\omega_0\in\partial\D$ such that
\begin{equation}\label{equ:limit:point}
    \sup_{\omega\in I(\omega_0,\delta)\setminus\{\omega_0\}}
    d_I(g(0),g(\omega))=\infty
    \qquad\text{for every } \delta\in(0,\pi).
\end{equation}
Indeed, let $M:=\max\{\diam(\Omega),1\}<\infty$. By the triangle inequality
$$d_I(x,y)\le  d_I(g(0),x)+d_I(g(0),y) \qquad\text{for all }x,y\in\Omega,$$
the assumption $\diam_I(\Omega)=\infty$ implies that the function $x\mapsto d_I(g(0),x)$ is unbounded. It is also continuous in $\Omega$ and attains the value $0$ at $x=g(0)$.
Since $\Omega$ is connected, the image of $\Omega$ under the function $x\mapsto d_I(g(0),x)$ is an interval containing $0$ and hence equals $[0,\infty)$. Therefore, for each $n\ge1$, we may choose
$x_n\in\Omega$ such that
\[
   d_I(g(0),x_n)=(3n+1)M.
\]
For each $x_n\in\Omega$, we choose a nearest boundary point $\eta_n\in\partial\Omega$:
\[
    |x_n-\eta_n|=\dist(x_n,\partial\Omega).
\]
Then $[x_n,\eta_n)\subset\Omega$, and hence
\[
    d_I(x_n,\eta_n)=|x_n-\eta_n|\le \diam(\Omega) \le M.
\]
The triangle inequality gives
\[
    3nM\le d_I(g(0),\eta_n)\le(3n+2)M<\infty.
\]
Since the intervals $[3nM,(3n+2)M]$, $n\ge1$, are pairwise disjoint,
the points $\eta_n$ are distinct.

Let $\omega_n:=g^{-1}(\eta_n)\in\partial\D$.
By compactness of $\partial\D$, a subsequence $\{\omega_{n_k}\}$ converges to
some $\omega_0\in\partial\D$.
Since the boundary points $\omega_n$ are distinct, we may discard
at most one term to ensure that $\omega_{n_k}\ne\omega_0$ for every $k$.
For each $\delta\in(0,\pi)$, we therefore have
$\omega_{n_k}\in I(\omega_0,\delta)\setminus\{\omega_0\}$
for all sufficiently large $k$, while
    $$d_I(g(0),g(\omega_{n_k}))
    =d_I(g(0),\eta_{n_k})
    \ge3n_kM\longrightarrow\infty.$$
This proves \eqref{equ:limit:point}.

Let $\omega_0$ be as in \eqref{equ:limit:point}. After rotating $\mathbb D$ and, if necessary, reflecting both
$\mathbb D$ and $\Omega$ across the real axis, we may assume that  $\omega_0=e^{i\pi}$ and
\begin{equation}\label{ineq: d_I to infty}
    \sup_{\pi-\delta<t<\pi}
    d_I(g(0),g(e^{it}))=\infty
    \qquad\text{for every } \delta\in(0,\pi).
\end{equation}

For every $\omega\in\partial\D$, $g([0,\omega))\subset\Omega$ is a hyperbolic geodesic joining $g(0)$ and $g(\omega)$. Applying Lemma~\ref{lem: gehring-hayman} and passing to the boundary
by a limiting argument, we obtain
\begin{equation}\label{ineq: the constant in GH lemma}
    \ell(g([0,\omega]))
    \le C\,d_I(g(0),g(\omega))
\end{equation}
where $C$ is an absolute constant.
Moreover, for any sequence $\omega_j\to\omega$ in $\partial\D$,
Fatou's lemma gives
\[
    \ell(g([0,\omega]))
    =\int_0^1 |g'(r\omega)|\,dr
    \le \liminf_{j\to\infty}
    \int_0^1 |g'(r\omega_j)|\,dr
    =\liminf_{j\to\infty}\ell(g([0,\omega_j])).
\]
Hence the function $\omega\mapsto\ell(g([0,\omega]))$
is lower semicontinuous on $\partial\D$.

Next, we construct by induction a sequence of intervals
   $\{I(\theta_n,\delta_n)\}_{n\ge 0}$
    with the following properties:
    \begin{itemize}
        \item[(1)] The intervals $\{I(\theta_n,\delta_n)\}_{n\in\mathbb N}$ are pairwise disjoint with $\omega_0\notin I(\theta_n,\delta_n)$;
        \item[(2)] For all $n\ge 0$, $\theta_n = e^{it_n}$, where $0 < t_1 < t_2 < \ldots < \pi$ and
        $$\lim_{n\to\infty} t_n = \pi;$$
        \item[(3)] For all $n \ge 0$, it holds that
        \begin{equation}\label{equ:interval:lowerbound}
            \inf_{\omega\in I(\theta_n,\delta_n)} d_I(g(0),g(\omega)) \geq 4^n.
        \end{equation}
    \end{itemize}
    Set $t_0=\delta_0=0$ and suppose that the  $I(t_0,\delta_0),\cdots,I(t_{n-1},\delta_{n-1})$ have been chosen. By \eqref{ineq: d_I to infty}, we may choose
$t_{n-1}+\delta_{n-1}<t_n<\pi$ such that
\[
    \pi-t_n<\frac{\pi}{2^n}
    \quad\text{and}\quad
    \ell(g([0,e^{it_n}]))>2C4^n,
\]
where $C$ is the absolute constant in \eqref{ineq: the constant in GH lemma}.
Hence, $t_n$ is increasing and satisfies $\lim_{n\to\infty}t_n=\pi$. Set $\theta_n=e^{it_n}$.
By lower semicontinuity, we may choose $\delta_n>0$ small enough such that
\[
    t_{n-1}+\delta_{n-1}
    <t_n-\delta_n<t_n+\delta_n<\pi
    \quad\text{and}\quad\ell(g([0,\omega]))>C4^n
\]
for every $\omega\in I(\theta_n,\delta_n)$. Combining this with \eqref{ineq: the constant in GH lemma}, we have 
\[
   \omega_0=e^{i\pi}\notin I(\theta_n,\delta_n)\quad\text{and}\quad
   d_I(g(0),g(\omega)) \geq 4^n
\]
for every $\omega\in I(\theta_n,\delta_n)$, which proves \eqref{equ:interval:lowerbound} and completes the construction of $I(\theta_n,\delta_n)$.

We now construct the boundary homeomorphism $\varphi$. For each $n\ge1$, we set
\[
    A_n:=\left\{e^{it}:
        \pi-\frac{\pi}{2^n}
        \le t\le
        \pi-\frac{\pi}{2^n}+\frac{\pi}{2\cdot4^n}
    \right\}.
\]
and define the boundary homeomorphism $h$ as follows: $h(e^{i\pi})=e^{i\pi}$ and $h$ maps each $A_n$ onto $I(\theta_n,\delta_n)$ by an
orientation-preserving affine map in angular coordinates. We then extend $h$ in the same way across the complementary arcs.
Since $A_n$ and $I(\theta_n,\delta_n)$
have the same circular order and accumulate only at $e^{i\pi}$,
$h$ is indeed a sense-preserving homeomorphism, and the desired boundary homeomorphism $\varphi$ is obtained by setting  $\varphi:=g\circ h$.
\begin{figure}[htbp]
    \centering
    \scalebox{0.8}
    {%
\usetikzlibrary{arrows.meta,decorations.markings}
\begin{tikzpicture}[
  font=\small,
  base/.style={draw=blue,line width=.5pt},
  special/.style={draw=red,line width=1.65pt,line cap=round},
  directed/.style={
    postaction={decorate},
    decoration={markings,mark=at position .7 with {
      \pgfsetstrokecolor{black}\pgfsetfillcolor{black}
      \arrow{Stealth[length=1.5mm,width=1.15mm]}
    }}
  },
  leader/.style={draw=black!80,line width=.4pt},
  arc label/.style={font=\footnotesize,inner sep=1.2pt},
  every node/.style={inner sep=2pt}
]
\def\R{2.25}
\def\Sep{7.5}

% Source circle.
\draw[base] (0,0) circle[radius=\R];
\foreach \n in {1,2,3,4,5,6}{
  \pgfmathsetmacro{\aa}{180-180/(2^\n)}
  \pgfmathsetmacro{\bb}{\aa+90/(4^\n)}
  \ifnum\n<3
    \draw[special,directed] (\aa:\R)
      arc[start angle=\aa,end angle=\bb,radius=\R];
  \else
    \draw[special] (\aa:\R)
      arc[start angle=\aa,end angle=\bb,radius=\R];
  \fi
  \ifnum\n<4
    \fill[red] (\aa:\R) circle[radius=.6pt];
    \fill[red] (\bb:\R) circle[radius=.6pt];
  \fi
}
\draw[leader] (101.25:\R) -- (101.25:2.52)
  node[above,arc label] {$A_1$};
\draw[leader] (135.1:\R) -- (135.1:2.52)
  node[left,arc label] {$A_2$};
\draw[leader] (158.203125:\R) -- (158.203125:2.52)
  node[left,arc label] {$A_3$};
\fill (-\R,0) circle[radius=1.2pt];
\node[below left=3pt] at (-\R,0) {$e^{i\pi}$};
\node at (-2.53,.4) {$\vdots$};
\draw[base,directed]
  (292:\R) arc[start angle=292,end angle=316,radius=\R];
\node[below] at (0,-2.53) {$\partial\mathbb D$};

% Target circle: these example arcs have the same circular order.
\begin{scope}[xshift=\Sep cm]
  \draw[base] (0,0) circle[radius=\R];
  \foreach \n/\aa/\bb in {
      1/92/124,2/142/156,3/166/172,
      4/175/176.5,5/178/178.6,6/179.1/179.35}{
    \ifnum\n<3
      \draw[special,directed] (\aa:\R)
        arc[start angle=\aa,end angle=\bb,radius=\R];
    \else
      \draw[special] (\aa:\R)
        arc[start angle=\aa,end angle=\bb,radius=\R];
    \fi
    \ifnum\n<4
      \fill[red] (\aa:\R) circle[radius=.6pt];
      \fill[red] (\bb:\R) circle[radius=.6pt];
    \fi
  }
  \draw[leader] (108:\R) -- (108:2.52)
    node[above,arc label] {$I(\theta_1,\delta_1)$};
  \draw[leader] (146:\R) -- (146:2.52)
    node[left,arc label] {$I(\theta_2,\delta_2)$};
  \draw[leader] (169:\R) -- (-2.45,.78)
    node[left,arc label] {$I(\theta_3,\delta_3)$};
  \fill (-\R,0) circle[radius=1.2pt];
  \node[below left=3pt] at (-\R,0) {$e^{i\pi}$};
  \node at (-2.53,.4) {$\vdots$};
  \draw[base,directed]
    (292:\R) arc[start angle=292,end angle=316,radius=\R];
  \node[below] at (0,-2.53) {$\partial\mathbb D$};
\end{scope}

% The map and its fixed accumulation point.
\draw[draw=black,-{Stealth[length=2.5mm,width=1.8mm]},line width=.7pt]
  (2.7,0) -- (4.8,0) node[midway,above=4pt] {$h$};
\node at (3.75,-1.0) {$h(e^{i\pi})=e^{i\pi}$};
\node at (3.75,-1.43) {$h(A_n)=I(\theta_n,\delta_n)$};
\end{tikzpicture}
    }
    \caption{The construction of $h$}
    \label{fig: The construction of $g$}
\end{figure}
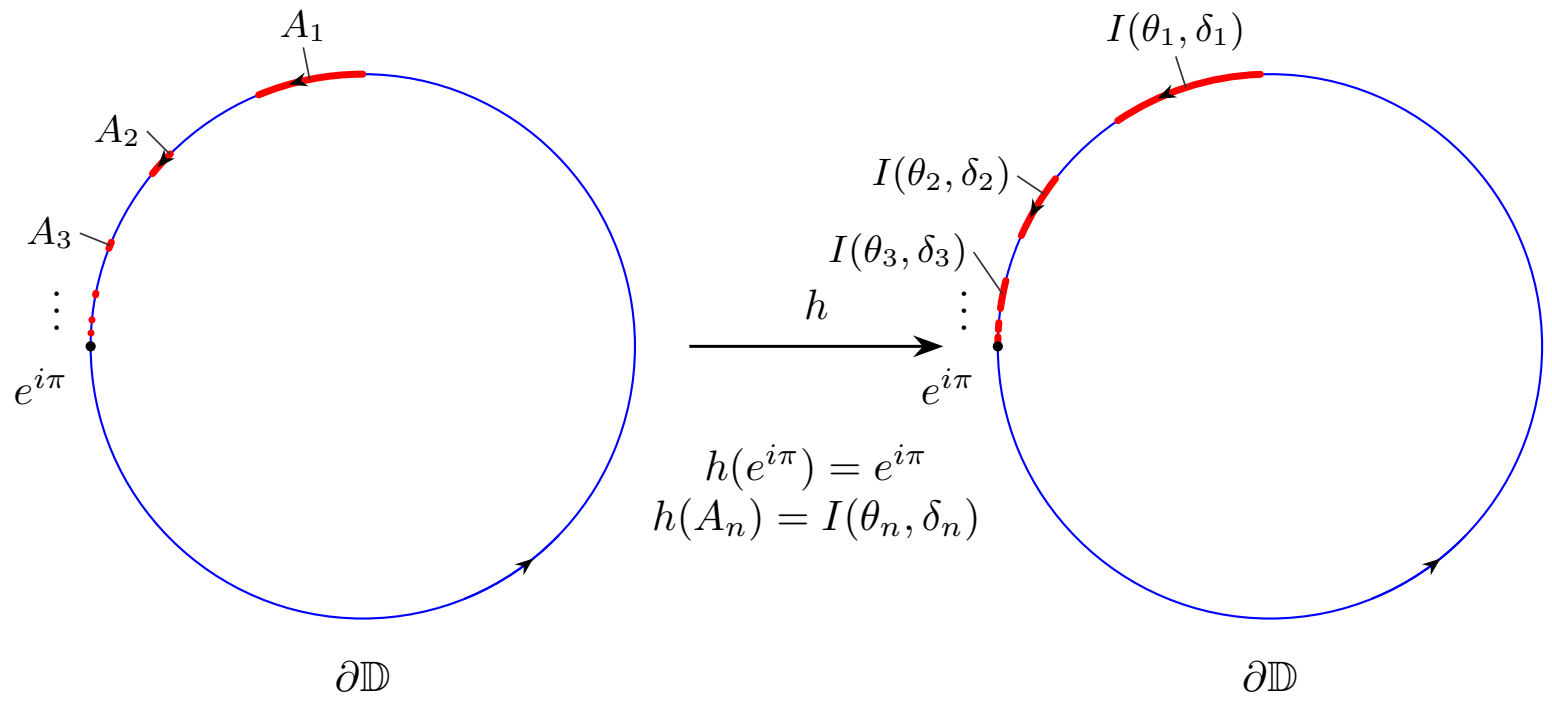

We finally show that such a $\varphi$ does not admit a homeomorphic-$BV(\D;\C)$ extension. 
Suppose, to the contrary, that such an extension $f_\varphi$ exists. Since $f_\varphi(\mathbb S(0,\frac{1}{2}))$ is compact in $\Omega$ and
$d_I(g(0),\cdot)$ is continuous in $\Omega$,
the internal distance $d_I(g(0),f_\varphi(\tfrac12\omega))$ is uniformly bounded for $\omega\in\partial\D$. 
Moreover, by the definition of $\varphi$ and \eqref{equ:interval:lowerbound} , we have
$h(A_n)=I(\theta_n,\delta_n)$ and 
\[
    d_I(g(0),\varphi(\omega))\ge4^n
\]
for every $\omega\in A_n$ and every $n\ge1$. Consequently, there exists an integer $n_0\ge1$ such that,
\[
    d_I(f_\varphi(\tfrac12\omega),\varphi(\omega))\ge d_I(g(0),\varphi(\omega))-d_I(g(0),f_\varphi(\tfrac12\omega))\ge \frac{4^n}{2}
\]
for every $\omega\in A_n$ and every $n\ge n_0$.

Applying the BV slicing theorem (Theorem \ref{thm:sobolev-slicing}) and using
the continuity of $f_\varphi$ on $\overline\D$,
the curve $\gamma_t(r)=f_\varphi(re^{it})$ has bounded variation
on $1/2<r<1$ for almost every $t$. In particular, 
its total variation equals its length and is therefore
at least the internal distance between its endpoints
\[
    \operatorname{Var}_{[\frac{1}{2},1]}(f_\varphi)=\ell(\gamma_t[{\scriptstyle\frac{1}{2}},1])\ge d_I\bigl(f_\varphi(\tfrac12e^{it}),\varphi(e^{it})\bigr)
\]
Consequently,
\begin{multline*}
    |Df_\varphi|(\D)
    \ge \frac{1}{2}\int_0^{2\pi}
        \operatorname{Var}_{[\frac{1}{2},1]}(f_\varphi)\,dt
    \ge \frac12\int_0^{2\pi}
        d_I\bigl(f_\varphi(\tfrac12e^{it}),\varphi(e^{it})\bigr)\,dt\\
    \ge \frac14\sum_{n=n_0}^{\infty}
        4^n\mathcal H^1(A_n)
     = \frac{\pi}{8}\sum_{n=n_0}^{\infty}1
     = \infty.
\end{multline*}
This contradicts $f_\varphi\in BV(\D;\C)$ and completes the proof.
\end{proof}

We are now ready to prove Theorem~\ref{thm:intro-all-boundary-data}.

\begin{proof}[Proof of Theorem~\ref{thm:intro-all-boundary-data}]
Suppose that $\diam_I(\Omega)<\infty$. Then every boundary homeomorphism
$\varphi\colon\partial\D\to\partial\Omega$ has a homeomorphic
$W^{1,1}$-extension by Proposition~\ref{thm: W11-extension}. Conversely, suppose that  $\diam_I(\Omega)=\infty$. By
Proposition~\ref{main-3}, there exists a boundary homeomorphism
$\varphi\colon\partial\D\to\partial\Omega$ that does not have a
homeomorphic $BV$-extension. Since
$W^{1,1}(\D;\C)\subset BV(\D;\C)$, this $\varphi$ cannot have a
homeomorphic $W^{1,1}$-extension. This proves the theorem.
\end{proof}

\subsection{Uniform bounds and equivalent extension properties} 
In this subsection, we establish Corollaries \ref{cor:uniform-W11-extension} and \ref{thm: equivalent-extension-properties}. More precisely, we give
equivalent characterizations of the property that every boundary
homeomorphism admits a homeomorphic $W^{1,1}$-extension.

\begin{proof}[Proof of Corollary \ref{cor:uniform-W11-extension}]
The implication $\textup{(ii)}\Rightarrow\textup{(i)}$
is immediate. Conversely, suppose that  \textup{(i)} holds.
Theorem~\ref{thm:intro-all-boundary-data} implies that $M_\Omega:=\operatorname{diam}_I(\Omega)<\infty$. By Proposition \ref{thm: W11-extension} and Remark \ref{rem: bi-Sobolev}, every boundary
homeomorphism $\varphi$
admits a homeomorphic $W^{1,1}$-extension
$f_\varphi:\overline{\mathbb D}\to\overline\Omega$
satisfying
\[
    \int_\Omega|Df^{-1}_\varphi(w)|\, dA(w)+\int_{\mathbb D}|Df_\varphi(z)|\,dA(z)
    \le C M_\Omega,
\]
where $C$ is an absolute constant.
Since $f_\varphi(\D)=\Omega$, we obtain
\[
    \|f_{\varphi}\|_{W^{1,1}(\mathbb D;\Omega)}+\|f^{-1}_{\varphi}\|_{W^{1,1}(\Omega;\D)}
    \le \pi\sup_{w\in\Omega}|w|+|\Omega|+C M_\Omega
    =:C_\Omega<\infty.
\]
The constant $C_\Omega$ depends only on $\Omega$,
which proves \textup{(ii)}.
\end{proof}

\begin{proof}[Proof of Corollary \ref{thm: equivalent-extension-properties}]
The implications
$\textup{(i)}\Rightarrow\textup{(iii)}$,
$\textup{(ii)}\Rightarrow\textup{(iii)}$, and
$\textup{(v)}\Rightarrow\textup{(iv)}$
follow directly from the definitions,
and the inclusion $W^{1,1}(\D;\Omega)\subset BV(\D;\Omega)$ yields
$\textup{(iii)}\Rightarrow\textup{(iv)}$ and
$\textup{(ii)}\Rightarrow\textup{(v)}$. Moreover, Proposition~\ref{main-3} and Remark~\ref{rem: bi-Sobolev}
yield $\textup{(iv)}\Rightarrow\textup{(vi)}$ and
$\textup{(vi)}\Rightarrow\textup{(ii)}$, respectively.
It therefore remains to prove $\textup{(vi)}\Rightarrow\textup{(i)}$.

% \medskip
% \noindent
% \textbf{Proof of $\textup{(iv)}\Rightarrow\textup{(vi)}$.}
% Suppose that \textup{(iv)} holds but
% $\operatorname{diam}_I(\Omega)=\infty$.
% By Theorem~\ref{main-3}, there exists a boundary homeomorphism
% $\varphi:\partial\D\to\partial\Omega$
% that admits no homeomorphic $BV$-extension.
% This contradicts \textup{(iv)} and proves \textup{(vi)}.

\medskip
\noindent
\textbf{Proof of \textup{(vi)} $\Rightarrow$ \textup{(i)}.}
Fix an orientation-preserving boundary homeomorphism
$\varphi\colon\partial\D\to\partial\Omega$, and define
$g$, $\psi$, $\Psi$, and $f_\varphi=g\circ\Psi$
as in the proof of Proposition \ref{thm: W11-extension}.
Then $f_\varphi\in W^{1,1}(\D;\Omega)$.
By Proposition~\ref{prop:circle-extension} and the conformality
of $g$, the map $f_\varphi$ is locally bi-Lipschitz.
Since it is orientation-preserving, we have
$J_{f_\varphi}>0$ a.e. in $\D$.
It follows from the area formula that
   $$ \int_{\D}J_{f_\varphi}(z)\,dA(z)\le|\Omega|<\infty.$$
In particular, we have
    $$K(z):=\frac{|Df_\varphi(z)|^2}{J_{f_\varphi}(z)}\in [1,\infty)\qquad\text{a.e. in }\D.$$
Defining $K=1$ on the exceptional null set, we conclude that
$f_\varphi$ is a homeomorphic extension of finite distortion.
This proves \textup{(i)}.

\end{proof}

\subsection{Sharpness of the \texorpdfstring{$W^{1,1}$}{W1,1}-extension}
In this subsection, we prove Corollary \ref{cor:intro-sharpness}, according to which the $W^{1,1}$-regularity in Theorem \ref{thm:intro-all-boundary-data} cannot be improved
to $W^{1,p}$-regularity for any $p>1$.

To prove Corollary \ref{cor:intro-sharpness}, we first construct
an example on the unit square.
\begin{prop}\label{prop: sharpness-W11}
Let $Q:=(0,1)^2$ be the unit square. There exist a Jordan domain $\Omega\subset\mathbb C$ with
$\diam_I(\Omega)<\infty$ and a boundary homeomorphism
$\varphi\colon\partial Q\to\partial\Omega$ which admits no
homeomorphic $W^{1,p}$-extension for any $p>1$.
\end{prop}
\begin{proof}
    We begin by constructing the target domain $\Omega$. For every $j\ge2$, we set
    $$a_j:=2^{-j},\quad
    b_j:=\tfrac{5}{4}\,2^{-j},\quad
    h_j:=(4\log j)^{-1},$$
and define
\[
    \Omega
    :=Q\setminus\bigcup_{j=2}^{\infty}
    \bigl([a_j,b_j]\times[0,h_j]\bigr).
\]
By construction, the boundary of $\Omega$ is a Jordan curve. Moreover, any two points of $\Omega$ can be joined  by at most two vertical segments and one horizontal segment. Hence, we have $\operatorname{diam}_I(\Omega)\le3$.

We choose $\delta>0$ sufficient small such that 
\[
    6\sum_{j=2}^{\infty}s_j<1/4,
    \quad\text{where}\quad s_j:=\frac{\delta}{j(\log j)^2},
\]
and set
\[
    c_j:=6\sum_{k>j}s_k+2s_j
    \quad\text{and}\quad
    I_j:=[c_j-2s_j,c_j+2s_j].
\]
Clearly, $I_j\subset[0,1/4]$ for all $j\ge2$. We then define a boundary homeomorphism
$\varphi\colon\partial Q\to\partial\Omega$ as follows. 
For each $j\ge2$, we map the subintervals $[c_j-2s_j,c_j-s_j]$ and $[c_j+s_j,c_j+2s_j]$ affinely onto the segments between
$$ (a_j,0),\,\,(a_j,h_j/2) \quad \text{and } \quad (b_j,h_j/2),\,\, (b_j,0),$$
respectively, and map the middle subinterval $[c_j-s_j,c_j+s_j]$ piecewise affinely onto the polygonal arc determined by 
$$(a_j,h_j/2),\quad (a_j,h_j),\quad
    (b_j,h_j),\quad (b_j,h_j/2);$$
see Figure \ref{fig: partial Q to partial Omgea}. We also map the remaining intervals on the bottom side of $Q$ affinely onto the corresponding bottom boundary segments of $\Omega$, and let $\varphi$ be the identity on the other three sides of $Q$. Since $a_j,b_j,h_j\to0$, the resulting map is continuous at the origin, and hence is a homeomorphism.

\begin{figure}[htbp]
    \centering
    \scalebox{0.8}{%
\usetikzlibrary{arrows.meta,decorations.pathreplacing}

\colorlet{leftpart}{blue}
\colorlet{middlepart}{red}
\colorlet{rightpart}{blue}

\begin{tikzpicture}[
  x=1cm,y=1cm,
  every node/.style={font=\small},
  arc/.style={line width=1pt},
  directed/.style={arc,-{Stealth[length=1.8mm]}},
  brace/.style={decorate,decoration={brace,amplitude=4pt}},
  point/.style={circle,fill=black,inner sep=1.25pt}
]
  % The source interval is divided in the ratio 1:2:1.
  \coordinate (A) at (0,6);
  \coordinate (B) at (2,6);
  \coordinate (C) at (6,6);
  \coordinate (D) at (8,6);
  \draw[gray] (-.8,6) -- (A) (D) -- (8.8,6);
  \draw[directed,leftpart] (A) -- (B);
  \draw[directed,middlepart] (B) -- (C);
  \draw[directed,rightpart] (C) -- (D);
  \foreach \p in {A,B,C,D} \node[point] at (\p) {};
  \node[below=6pt] at (A) {$c_j-2s_j$};
  \node[below=6pt] at (B) {$c_j-s_j$};
  \node[below=6pt] at (C) {$c_j+s_j$};
  \node[below=6pt] at (D) {$c_j+2s_j$};
  \draw[brace] (0,6.65) -- (8,6.65)
    node[midway,above=6pt] {$I_j\subset\partial Q$};
  \draw[brace,decoration={brace,mirror,amplitude=4pt}]
    (0,5.25) -- (2,5.25) node[midway,below=6pt] {$s_j$};
  \draw[brace,decoration={brace,mirror,amplitude=4pt}]
    (2,5.25) -- (6,5.25) node[midway,below=6pt] {$2s_j$};
  \draw[brace,decoration={brace,mirror,amplitude=4pt}]
    (6,5.25) -- (8,5.25) node[midway,below=6pt] {$s_j$};

  \draw[-{Stealth[length=2mm]}] (4,4.4) -- (4,3.65)
    node[midway,right=4pt] {$\varphi$};

  % A local schematic of the target boundary; not drawn to scale.
  \fill[gray!7] (-.8,0) rectangle (8.8,3.45);
  \fill[white] (2,0) rectangle (6,3);
  \coordinate (Ap) at (2,0);
  \coordinate (Bp) at (2,1.5);
  \coordinate (Cp) at (6,1.5);
  \coordinate (Dp) at (6,0);
  \draw[gray] (-.8,0) -- (Ap) (Dp) -- (8.8,0);
  \draw[directed,leftpart] (Ap) -- (Bp);
  \draw[directed,middlepart] (Bp) -- (2,3) -- (6,3) -- (Cp);
  \draw[directed,rightpart] (Cp) -- (Dp);
  \draw[directed,middlepart] (3.6,3) -- (4.4,3);
  \foreach \p in {Ap,Bp,Cp,Dp} \node[point] at (\p) {};
  \node[below left=3pt] at (Ap) {$(a_j,0)$};
  \node[left=5pt] at (Bp) {$(a_j,h_j/2)$};
  \node[right=5pt] at (Cp) {$(b_j,h_j/2)$};
  \node[below right=3pt] at (Dp) {$(b_j,0)$};
  \node[above left=3pt] at (2,3) {$(a_j,h_j)$};
  \node[above right=3pt] at (6,3) {$(b_j,h_j)$};
  \node at (7.9,2.55) {$\Omega$};
  \node[right] at (8.8,0) {$\partial\Omega$};
  \end{tikzpicture}
    }
    \caption{The restriction of $\varphi$ to $I_j$.}
    \label{fig: partial Q to partial Omgea}
\end{figure}
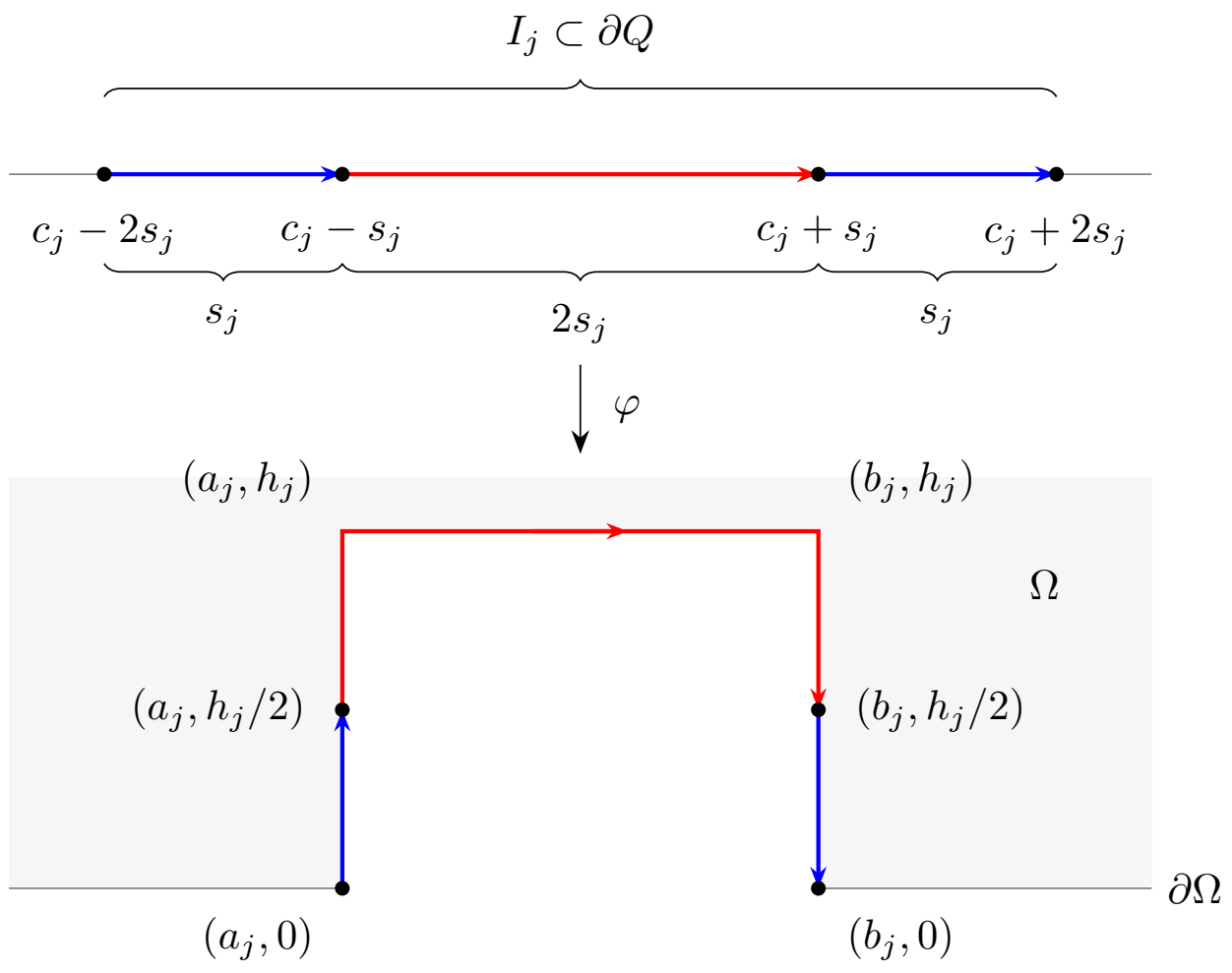

We now suppose that $\varphi$ admits a homeomorphic-$W^{1,p}$ extension $f_\varphi\colon\overline {Q}\to\overline{\Omega}$ for some $p>1$. For each $j\ge2$, we consider the
half-annulus
\[
    A_j
    :=\{c_j+re^{it}:s_j<r<2s_j,\ 0\le t\le\pi\}.
\]
and the semicircles
\[
    \gamma_{j,r}(t):=c_j+re^{it},\quad s_j<r<2s_j,\  0\le t\le\pi.
\]
By the Sobolev slicing theorem (Theorem \ref{thm:sobolev-slicing}), the curve $f_\varphi\circ\gamma_{j,r}$
is absolutely continuous for a.e. $s_j<r<2s_j$. Noting that the endpoints of $f_\varphi\circ\gamma_{j,r}$ lie on the vertical segments between
$$\quad (a_j,0), \,\,(a_j,h_j/2) \quad \text{and} \quad (b_j,h_j/2),\,\, (b_j,0),$$
respectively,
it follows that $h_j\le\ell(f_\varphi\circ\gamma_{j,r})$, see Figure \ref{fig: the image of f circ gamma}. 
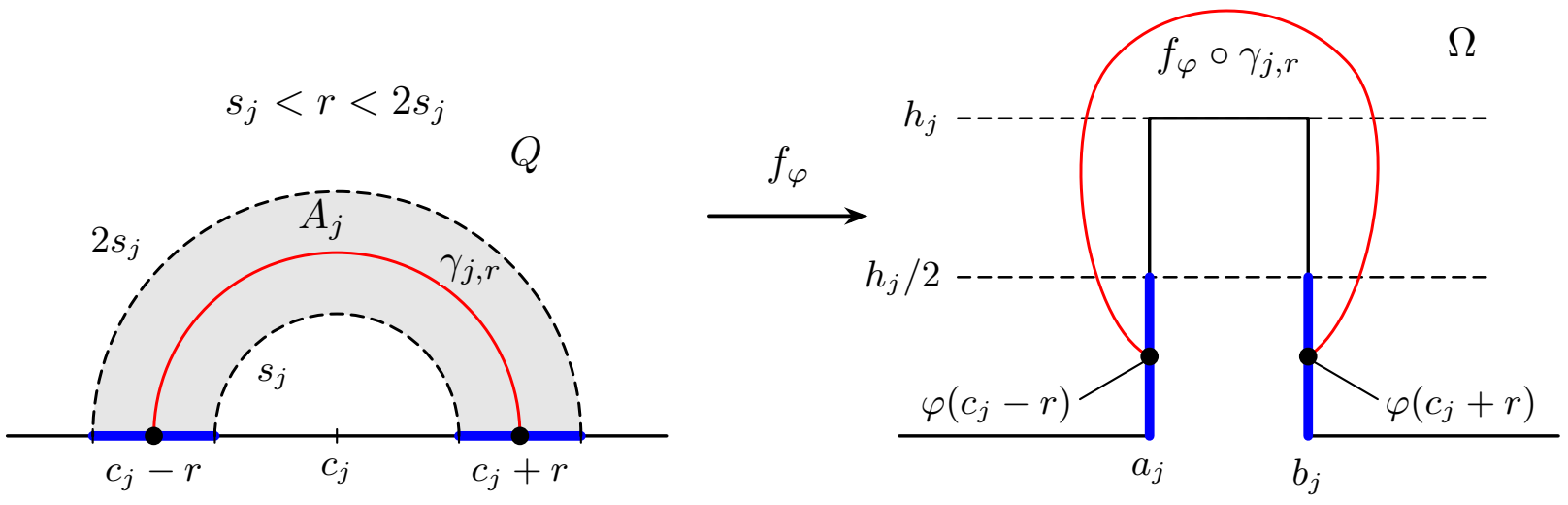
\begin{figure}[htbp]
    \centering
    \scalebox{0.8}
    {%
\begin{tikzpicture}[x=1.04cm,y=1.04cm,
  font=\small,line cap=round,line join=round,
  >={Stealth[length=2mm]}]
  % Source: a local view of the lower side of Q.
  \path[fill=annulusfill,even odd rule]
    (2,0) arc[start angle=0,end angle=180,radius=2] -- cycle
    (1,0) arc[start angle=0,end angle=180,radius=1] -- cycle;
  \draw[boundary,line width=.8pt] (-2.7,0)--(2.7,0);
  \draw[endpoint,line width=2.3pt] (-2,0)--(-1,0);
  \draw[endpoint,line width=2.3pt] (1,0)--(2,0);
  \draw[guide,densely dashed,line width=.7pt]
    (2,0) arc[start angle=0,end angle=180,radius=2];
  \draw[guide,densely dashed,line width=.7pt]
    (1,0) arc[start angle=0,end angle=180,radius=1];
  \draw[curve,line width=.7pt]
    (1.5,0) arc[start angle=0,end angle=180,radius=1.5];
  \node[boundary,fill=annulusfill,inner sep=0.0pt] at (1.1,1.36)
    {$\gamma_{j,r}$};
  \node[boundary] at (-.12,1.77) {$A_j$};
  \node[guide,fill=white,inner sep=1pt,font=\footnotesize]
    at (-1.8,1.57) {$2s_j$};
  \node[guide,fill=white,inner sep=1pt,font=\footnotesize]
    at (-.52,.48) {$s_j$};
  \node at (0,2.70) {$s_j<r<2s_j$};
  \node at (1.55,2.30) {$Q$};
  \foreach \x in {-2,-1,0,1,2}
    \draw[boundary,line width=.5pt] (\x,-.055)--(\x,.055);
  \node[below,font=\footnotesize] at (0,-.06) {$c_j$};
  \fill[boundary] (-1.5,0) circle[radius=2.2pt];
  \fill[boundary] (1.5,0) circle[radius=2.2pt];
  \node[boundary,below,font=\footnotesize] at (-1.5,-.06) {$c_j-r$};
  \node[boundary,below,font=\footnotesize] at (1.5,-.06) {$c_j+r$};
  % Mapping arrow.
  \draw[->,boundary,line width=.9pt] (3.05,1.8)--(4.35,1.8)
    node[midway,above=3pt] {$f_\varphi$};
  % Target: a local view of the j-th removed rectangle.
  \begin{scope}[xshift=7.6cm]
    \draw[guide,densely dashed,line width=.6pt] (-2.22,1.3)--(2.13,1.3);
    \draw[guide,densely dashed,line width=.6pt] (-2.22,2.6)--(2.13,2.6);
    \node[guide,left,font=\footnotesize] at (-2.23,1.3) {$h_j/2$};
    \node[guide,left,font=\footnotesize] at (-2.23,2.6) {$h_j$};
    \draw[boundary,line width=.8pt]
      (-2.7,0)--(-.65,0)--(-.65,2.6)--(.65,2.6)--(.65,0)--(2.7,0);
    \draw[endpoint,line width=2.3pt] (-.65,0)--(-.65,1.3);
    \draw[endpoint,line width=2.3pt] (.65,0)--(.65,1.3);
    \node at (1.91,3.22) {$\Omega$};
    \draw[curve,line width=.7pt]
      (-.65,.65)
      .. controls (-1.19,.93) and (-1.45,2.58) .. (-.93,3.10)
      .. controls (-.44,3.62) and (.46,3.61) .. (.98,3.05)
      .. controls (1.46,2.51) and (1.16,.99) .. (.65,.65);
    \node[boundary,fill=white,inner sep=1.1pt] at (0,3.1)
      {$f_\varphi\circ\gamma_{j,r}$};
    \fill[boundary] (-.65,.65) circle[radius=2.2pt];
    \fill[boundary] (.65,.65) circle[radius=2.2pt];
    \draw[boundary,line width=.5pt] (-.69,.61)--(-1.22,.30);
    \draw[boundary,line width=.5pt] (.69,.61)--(1.22,.30);
    \node[boundary,anchor=east,font=\footnotesize] at (-1.15,.26)
      {$\varphi(c_j-r)$};
    \node[boundary,anchor=west,font=\footnotesize] at (1.15,.26)
      {$\varphi(c_j+r)$};
    \node[below,font=\footnotesize] at (-.65,-.06) {$a_j$};
    \node[below,font=\footnotesize] at (.65,-.06) {$b_j$};
  \end{scope}
  \end{tikzpicture}}
  \caption{The image of $f_\varphi\circ\gamma_{j,r}$.}
   \label{fig: the image of f circ gamma}
\end{figure}
Consequently, by applying H\"older's inequality, we obtain
\[
    h_j
    \le \ell(f_\varphi\circ\gamma_{j,r})
    \le \int_{\gamma_{j,r}}|Df_\varphi|\,ds
    \le (\pi r)^{1-1/p}
       \left(\int_{\gamma_{j,r}}|Df_\varphi|^p\,ds\right)^{1/p}, 
\]
where we used $\ell(\gamma_{r,j})=\pi r$ in the last inequality. Integrating both sides with respect to $r$ from $s_j$ to $2s_j$ and using $(\pi r)^{1-p}\ge(2\pi s_j)^{1-p}$, we obtain
\[
    \int_{A_j}|Df_\varphi(z)|^p\,dA(z)
    \ge (2\pi)^{1-p}h_j^p s_j^{2-p}.
\]
Since the half-annuli $A_j$ are pairwise disjoint, it follows that
\[
    \int_Q|Df_\varphi(z)|^p\,dA(z)
    \ge C_{p,\delta}
       \sum_{j=2}^{\infty}
       \frac{j^{p-2}}{(\log j)^{4-p}}
    =\infty
\]
for every $p>1$. This contradicts $f_\varphi\in W^{1,p}(Q,\Omega)$.
\end{proof}

\begin{proof}[Proof of Corollary \ref{cor:intro-sharpness}]
    Let the Jordan domain $\Omega$ and the homeomorphism $\widetilde{\varphi}\colon\partial Q\to\partial\Omega$
be those constructed in Proposition \ref{prop: sharpness-W11}. In particular, $\operatorname{diam}_I(\Omega)<\infty$. 
By Theorem \ref{thm:intro-all-boundary-data}, there exists a constant $C_\Omega<\infty$ only depending on $\Omega$ such that 
\[
    \|f_\varphi\|_{W^{1,1}(\D;\Omega)}\le C_\Omega.
\]

We then choose a bi-Lipschitz
homeomorphism $H\colon\overline\D\to\overline Q$ and define the homeomorphism $\varphi\colon\partial \D\to\partial\Omega$ by 
\[
    \varphi(z):=\widetilde{\varphi}(H(z)).
\]
Suppose that $\varphi$ admits a homeomorphic $W^{1,p}$-extension
for some $p>1$. Then $\widetilde{f_\varphi}:=f_\varphi\circ H^{-1}$
is a homeomorphic extension of $\widetilde{\varphi}$.
Since $H^{-1}$ is bi-Lipschitz, we have
$\widetilde{f_\varphi}\in W^{1,p}(Q,\Omega)$, which contradicts Proposition \ref{prop: sharpness-W11}.
\end{proof}

\section{The case of a single parametrization}\label{thm1.1}
In this section, we prove Theorem \ref{thm:intro-criterion}.

\begin{proof}[Proof of Theorem \ref{thm:intro-criterion}]

The implication \textup{(ii)}$\Rightarrow$\textup{(i)} is immediate.
It remains to prove $\textup{(i)}\Rightarrow\textup{(iii)}$ and
$\textup{(iii)}\Rightarrow\textup{(ii)}$.

\medskip
\noindent
\textbf{Proof of \textup{(i)} $\Rightarrow$ \textup{(iii)}.}
Suppose that $\varphi:\partial\D\to\partial\Omega$ admits a homeomorphic
$W^{1,1}$-extension
$f_\varphi:\cl{\D}\to\cl{\Omega}$, and fix $x_0\in\Omega$. Since the function $\dI(z,x_0)$ is continuous in $\Omega$ and $K:=f_{\varphi}(\mathbb S(0,1/2))\subset\Omega$ is compact, we have
$M:=\sup_{z\in K}\dI(z,x_0)<\infty$.

By the continuity of $f_{\varphi}$ up to $\partial\D$ and the Sobolev slicing theorem (Theorem \ref{thm:sobolev-slicing}), we know that for almost every $\xi\in\partial \D$, the curve $\gamma_\xi(r):=f_{\varphi}(r\xi)$ is absolutely continuous on ${\scriptstyle\frac{1}{2}}\le r\le 1$. Therefore, we deduce that
\begin{equation*}
    d_I(\varphi(\xi),x_0)\le d_I(f_{\varphi}(\xi/2),x_0)+d_I(f_{\varphi}(\xi),f_{\varphi}(\xi/2))\le M+\int_{{\scriptstyle\frac{1}{2}}}^1|Df_{\varphi}(r\xi)|dr,
\end{equation*}
By integrating this inequality over  $\partial\D$,  we conclude that
$$\int_{\partial \D}\dI\bigl(\varphi(\xi),x_0\bigr)\,d\mathcal{H}^1(\xi)\le 2\pi M+2\int_{\D}|Df_{\varphi}(z)|\,dA(z)<\infty.$$
Since $x_0$ is arbitrary, this proves
\textup{(i)}$\Rightarrow$\textup{(iii)}.

\medskip
\noindent
\textbf{Proof of \textup{(iii)} $\Rightarrow$ \textup{(ii)}.}
Suppose now that \textup{(iii)} holds for some $x_0\in\Omega$.
By precomposing $\varphi$ with the reflection $z\mapsto\overline z$
if necessary, we may assume that $\varphi$ is orientation-preserving.
Choose a conformal map $g\colon\D\to\Omega$ with $g(0)=x_0$.
By Carath\'eodory's theorem, $g$ extends to a homeomorphism
$g:\cl{\D}\to\cl{\Omega}$, which we denote by the same notation.

For every $\xi\in\partial\D$, we have
\begin{equation}\label{eq:boundary-gh}
\int_0^1|g'(r\xi)|\,dr
\le C\dI\bigl(x_0,g(\xi)\bigr).
\end{equation}
Indeed, suppose that $\dI\bigl(x_0,g(\xi)\bigr)<\infty$.
Let $\gamma\colon[0,1]\to\overline{\Omega}$ be a rectifiable curve
joining $x_0$ to $g(\xi)$ with $\gamma([0,1))\subset\Omega$, and choose $0<t_j<1$ with $t_j\to1$. Then we have the points
$w_j:=g^{-1}(\gamma(t_j))\in \D$ and $w_j\to\xi$.
By Lemma \ref{lem:conformal-invariance}, $g([0,w_j])$ is a hyperbolic geodesic in $\Omega$, and the Gehring--Hayman inequality (Lemma \ref{lem: gehring-hayman}) yields that 
$$\ell\bigl(g([0,w_j])\bigr)
\le C\ell\bigl(\gamma|_{[0,t_j]}\bigr)
\le C\ell(\gamma),$$
where $C$ is an absolute constant.
For $0<R<1$, since $w_j\to\xi\in\partial\D$, we have $|w_j|>R$ for all sufficiently large $j$, and hence
$$\int_0^R\left|g'\left(r\xi_j\right)\right|\,dr \le \int_0^{|w_j|}\left|g'\left(r\xi_j\right)\right|\,dr=\ell(g([0,w_j]))
\le C\ell(\gamma),$$
where $\xi_j:=w_j/|w_j|$.
By the uniform continuity of $g'$ on $\overline{B(0,R)}$,
we may let $j\to\infty$. Letting $R\to 1$ then yields
\[
    \int_0^1\left|g'\left(r\xi\right)\right|\,dr\le C\ell(\gamma)
\]
for all $\gamma\subset\Omega$ joining $x_0$ to $g(\xi)$, which proves \eqref{eq:boundary-gh}.

We now construct an extension of $\varphi$. Set
$\psi:=(g|_{\partial\D})^{-1}\circ\varphi\colon\partial\D\to\partial\D$,
and let $\Psi\colon\overline{\D}\to\overline{\D}$ be the homeomorphic
extension of $\psi$ constructed in Proposition \ref{prop:circle-extension}.
Then $f_{\varphi}:=g\circ\Psi\colon\overline{\D}\to\overline{\Omega}$
is a homeomorphic extension of $\varphi$ and it is locally Lipschitz in $\D$.

It remains to prove that $f_{\varphi}\in W^{1,1}(\D;\Omega)$.
Since $\Omega$ is bounded, it suffices to show that $|Df_{\varphi}|\in L^1(\D)$. By repeating the estimates leading to \eqref{eq:radial-energy} and
\eqref{eq:angular-energy} in the proof of
Proposition~\ref{thm: W11-extension},
we have 
\begin{multline*}
    \int_{\D} |Df_{\varphi}(z)|\,dA(z)
    \le \int_0^1\int_0^{2\pi}|\partial_rf_{\varphi}(re^{it})|\,r\,dt\,dr+ \int_0^1\int_0^{2\pi}|\partial_tf_{\varphi}(re^{it})|\,dt\,dr\\
    \le C\int_0^{1}\int_0^{2\pi}
|g'(r\psi(e^{it}))|\,dr\,dt
+C\int_0^1\int_0^{2\pi}|g'(re^{it})|\,r\,dt\,dr.
\end{multline*}

To estimate the first term, we apply Tonelli's theorem and
\eqref{eq:boundary-gh} with $\xi=\psi(e^{it})$ to obtain
\begin{equation}\label{ineq: g' first term}
\begin{aligned}
    \int_0^1\!\int_0^{2\pi}
    \bigl|g'(r\psi(e^{it}))\bigr|\,dt\,dr
    &\le C\int_0^{2\pi}
    \dI\bigl(x_0,g(\psi(e^{it}))\bigr)\,dt\\
    &=C\int_0^{2\pi}
    \dI\bigl(x_0,\varphi(e^{it})\bigr)\,dt
    \\
    &= C\int_{\partial\D}
    \dI\bigl(x_0,\varphi(\eta)\bigr)\,
    d\mathcal{H}^1(\eta)
    <\infty.
\end{aligned}
\end{equation}

For the second term, using polar coordinates, the
Cauchy--Schwarz inequality and the area formula,
we obtain
\begin{equation}\label{ineq: g' second term}
\begin{aligned}
    \int_0^1\!\int_0^{2\pi}
    |g'(re^{it})|\,r\,dt\,dr
    &= \int_\D |g'(z)|\,dA(z) \\
    &\le \sqrt{\pi}
    \left(\int_\D |g'(z)|^2\,dA(z)\right)^{1/2} \\
    &= \sqrt{\pi|\Omega|}
    < \infty.
\end{aligned}
\end{equation}
Combining \eqref{ineq: g' first term} and
\eqref{ineq: g' second term}, we conclude that
$f_\varphi\in W^{1,1}(\D;\Omega)$.

Finally, Proposition~\ref{prop:circle-extension} and the
conformality of $g$ show that $f_\varphi=g\circ\Psi$ is
locally bi-Lipschitz. The change of variables formula gives
$$\begin{aligned}
    \int_\Omega|D(f_\varphi^{-1})(w)|\,dA(w)
    &=\int_\D
    |(Df_\varphi(z))^{-1}|\,|J_{f_\varphi}(z)|\,dA(z)\\
    &=\int_\D|Df_\varphi(z)|\,dA(z)<\infty.
\end{aligned}$$
Since $f_\varphi^{-1}$ is bounded and $|\Omega|<\infty$, we conclude
that $f_\varphi^{-1}\in W^{1,1}(\Omega;\D)$.
Thus \textup{(iii)} implies \textup{(ii)}, completing the proof.
\end{proof}

\medskip
\textbf{Declaration of generative AI use:} During the preparation of this manuscript, the authors used ChatGPT to assist with proofreading and preparing the figures.

\bibliographystyle{amsalpha}
\bibliography{Sobolev-Homeo-exten.bib}
\end{document}